\documentclass[a4paper,11pt,oneside]{article}

\usepackage{fullpage}

\usepackage{amsmath}
\usepackage{amssymb,amsfonts}
\usepackage{amsthm}

\usepackage{graphicx}
\usepackage{xcolor}
\definecolor{darkred}{RGB}{153,0,0}
\definecolor{darkblue}{RGB}{0,0,153}
\definecolor{darkgreen}{RGB}{0,153,0}
\definecolor{forestgreen}{RGB}{0,110,51}
\definecolor{cb2red}{RGB}{228,26,28}
\definecolor{cb2blue}{RGB}{55,126,184}
\definecolor{cb2green}{RGB}{77,175,74}
\definecolor{cb2purple}{RGB}{152,78,163}
\definecolor{cb2orange}{RGB}{255,127,0}

\usepackage[colorlinks,
	citecolor=cb2blue,
	linkcolor=cb2blue,
	urlcolor=cb2blue]{hyperref}
\usepackage[ruled,vlined,linesnumbered,algosection,algo2e]{algorithm2e}
\usepackage[capitalize,noabbrev,nameinlink]{cleveref}
\crefname{line}{Step}{Steps}

\theoremstyle{plain}
\newtheorem{theorem}{Theorem}[section]

\theoremstyle{definition}
\newtheorem{assumption}[theorem]{Assumption}
\newtheorem{definition}[theorem]{Definition}
\theoremstyle{remark}
\newtheorem{remark}[theorem]{Remark}

\crefname{assumption}{Assumption}{Assumptions}

\newcommand{\KeywordsAnd}{\and}
\newcommand{\KeywordsEnd}{}
\newcommand{\keywords}[1]{\par\noindent{\small\def\and{\unskip,\ }{\bf Keywords. }#1.}\par}
\newcommand{\subclass}[1]{\par\noindent{\small\def\and{\unskip,\ }{\bf AMS MSC. }#1.}\par}

\newcommand\nnfootnote[1]{%
	\renewcommand\thefootnote{}\footnote{#1}%
	\addtocounter{footnote}{-1}%
}

\SetKw{KwOr}{or}
\SetKw{KwElse}{else}
\SetKw{KwAnd}{and}
\SetKw{KwReturn}{return}
\SetKw{KwGoTo}{go to}

\newcommand{\emailLink}[1]{\textsc{email} \href{mailto:#1}{#1}}
\newcommand{\orcidLink}[1]{\textsc{orcid} \href{https://orcid.org/#1}{#1}}
\newcommand{\amsmscLink}[1]{\href{http://www.ams.org/mathscinet/msc/msc2020.html?t=#1}{#1}}

\newcommand{\Z}{\mathbb{Z}}
\newcommand{\R}{\mathbb{R}}

\DeclareMathOperator{\proj}{proj}

\DeclareMathOperator{\dist}{dist}

\newcommand{\coloneqq}{:=}
\DeclareMathOperator*{\minimize}{minimize}

\DeclareMathOperator{\stt}{subject~to}
\DeclareMathOperator{\wrt}{over}

\newcommand{\func}[3]{#1 \colon #2 \to #3}
\newcommand{\innprod}[2]{\langle #1, #2 \rangle}
\newcommand{\normalcone}{\mathcal{N}}
\newcommand{\limnormalcone}{\normalcone^{\textup{lim}}}
\newcommand{\emphdef}[1]{\textsc{#1}}

\usepackage{tcolorbox}
\newtcolorbox{mybox}[1][]{%
	left=0pt,
	right=0pt,
	top=0pt,
	bottom=0pt,
	colback=gray!12,
	colframe=gray!12,
	width=\dimexpr\textwidth\relax,
	enlarge left by=0mm,
	boxsep=5pt,
	arc=1pt,outer arc=1pt,
	#1
}

\newcommand{\norm}[1]{\| #1 \|}
\newcommand{\closedball}{\mathbb{B}}
\newcommand{\lpfooter}{\textup{PL}}
\newcommand{\normlp}[1]{\norm{ #1 }_{\lpfooter}}
\newcommand{\lpball}{\closedball_{\lpfooter}}

\newcommand{\spaceX}{\R^n}
\newcommand{\spaceY}{\R^m}
\newcommand{\psimeas}{\Psi}
\newcommand{\lf}{\varphi}

\newcommand{\Lagr}{\mathcal{L}}

\newcommand{\XX}{\mathcal{X}}
\newcommand{\CC}{\mathcal{C}}
\newcommand{\Ybounded}{\mathcal{Y}}
\newcommand{\II}{\mathcal{I}}
\DeclareMathOperator{\totalvariation}{TV}

\newcommand{\TheAuthorADM}{Alberto De~Marchi}
\newcommand{\TheEmailADM}{alberto.demarchi@unibw.de}
\newcommand{\TheOrcidADM}{0000-0002-3545-6898}
\newcommand{\TheAuthorVN}{Viktoriya Nikitina}
\newcommand{\TheEmailVN}{viktoriya.nikitina@unibw.de}
\newcommand{\TheOrcidVN}{0000-0002-8523-7378}
\newcommand{\TheAuthorMG}{Matthias Gerdts}
\newcommand{\TheEmailMG}{matthias.gerdts@unibw.de}
\newcommand{\TheOrcidMG}{0000-0001-8674-5764}
\newcommand{\TheAffiliation}{%
	University of the Bundeswehr Munich,
	Department of Aerospace Engineering,
	Institute of Applied Mathematics and Scientific Computing,
	85577 Neubiberg, Germany%
}
\newcommand{\TheFunding}{This research was partially funded by dtec.bw -- Digitalization and Technology Research Center of the Bundeswehr [EMERGENCY-VRD, MissionLab, MORE, SeRANIS]. dtec.bw is funded by the European Union -- NextGenerationEU}
\newcommand{\ThePublicationNote}{This work has been accepted for presentation at the 13th IFAC Symposium on Nonlinear Control Systems (NOLCOS 2025) and for publication in the IFAC-PapersOnLine proceedings}
\newcommand{\TheTitle}{Hybrid optimal control with mixed-integer Lagrangian methods}
\newcommand{\TheKeywords}{%
	Hybrid dynamics\KeywordsAnd%
	Optimal control\KeywordsAnd
	Mixed-integer programming\KeywordsAnd%
	Constrained optimization\KeywordsAnd%
	Augmented Lagrangian\KeywordsEnd%
}
\newcommand{\TheAMSsubj}{%
	\amsmscLink{65K05}\and% Numerical mathematical programming methods
	\amsmscLink{90C06}\and% Large-scale problems in mathematical programming
	\amsmscLink{90C11}\and% Mixed integer programming
	\amsmscLink{90C30}% Nonlinear programming
}
\newcommand{\TheAbstract}{%
	Models involving hybrid systems are versatile in their application but difficult to optimize efficiently due to their combinatorial nature.
	This work presents a method to cope with hybrid optimal control problems which,
	in contrast to decomposition techniques, does not require relaxing the integrality constraints.
	Based on the discretize-then-optimize approach, our scheme addresses mixed-integer nonlinear problems under mild assumptions.
	The proposed numerical algorithm builds upon the augmented Lagrangian framework,
	whose subproblems are handled using successive mixed-integer linearizations with trust regions.
	We validate the performance of the numerical routine with extensive investigations using hybrid optimal control problems from different fields of application.
	Promising preliminary results are presented for
	a motion planning task with hysteresis and
	a Lotka-Volterra fishing problem with total variation.%
}%
\newcommand{\TheAcknowledgements}{%
	We thank Rebecca Richter for providing the dynamic programming model and implementation.%
}%

\begin{document}

\title{\bfseries \TheTitle}
\author{%
	\TheAuthorVN\thanks{%
		\emailLink{\TheEmailVN},
		\orcidLink{\TheOrcidVN}.%
	}\and%
	\TheAuthorADM\thanks{%
		\emailLink{\TheEmailADM},
		\orcidLink{\TheOrcidADM},
		corresponding author.%
	}\and%
	\TheAuthorMG\thanks{%
		\emailLink{\TheEmailMG},
		\orcidLink{\TheOrcidMG}.%
	}%
}%
\date{}

\maketitle
\begin{abstract}
	\TheAbstract
\end{abstract}
\subclass{\TheAMSsubj}

\section{Introduction}
\nnfootnote{The authors are with the \TheAffiliation.}
\nnfootnote{\TheFunding.}
\nnfootnote{\ThePublicationNote.}

We consider hybrid optimal control problems (OCPs), that is,
optimization problems involving possibly
nonsmooth dynamics, mixed-integer states or controls, or logical constraints.
Hybrid OCPs were investigated, e.g., by \cite{oldenburg2008disjunctive,grossmann2013systematic}
where hybrid dynamics are modeled with disjunctive programming.
An important subclass of hybrid OCPs are mixed-integer OCPs with discrete-valued controls.
Numerical methods often use variable time transformations, relaxation, or decomposition techniques; compare \cite{gerdts2006variable,sager2011combinatorial,buerger2020pycombina,plate2024second}
for important developments.
Such approaches are often tailored to discrete-valued controls and exploit the time structure.

In this work, we permit more general hybrid OCPs and apply a direct discretization method to transform the dynamic optimization problem into a finite dimensional mixed-integer nonlinear program (MINLP).
Our focus is not on the direct discretization method itself,
which exists in various versions,
compare \cite{gerdts2023optimal} for an overview,
nor on taking advantage of the time structure.
We investigate instead the numerical viability of an optimization scheme
that handles hybrid features as they stand,
namely without relaxing integrality constraints.
Moreover, since solving these problems to global optimality is in general computationally intractable, we do not seek global solutions and are content with suboptimal, yet practical ones, as in \cite{exler2007trust,quirynen2021sequential}.

We study a method for the numerical solution of a class of MINLPs of the following type:
\begin{align}
	\tag{P}\label{eq:P}
	\minimize~&f(x) &
	\stt~&x\in\XX ,\quad c(x) \in \CC
\end{align}
with $\func{f}{\XX}{\R}$ and $\func{c}{\XX}{\spaceY}$ differentiable functions,
$\CC\subset\spaceY$ a nonempty closed \emph{convex} set (projection-friendly in practice),
and $\XX$ a nonempty closed set with a mixed-integer linear (MIL) structure.
In particular, we assume $\XX$ to be described by linear inequalities and integrality
constraints on some variables, that is, 
\begin{equation} \label{EQ:MILS}
	\XX
	\coloneqq
	\left\{
	x \in \spaceX
	\,\middle|\,
	\begin{array}{l}
		A x \leq b ,\\
		x_i \in \Z \quad \forall i\in\II\subseteq\{1,\ldots,n\}
	\end{array}
	\right\}
\end{equation}
for some matrix $A$, vector $b$ and index set $\II$.
The general template \eqref{eq:P} covers hybrid OCPs
discretized with arbitrary schemes,
and can capture various features such as mixed-integer control inputs, nonsmooth dynamics and state jumps.

In the following, we will partition the decision variables $x$ into a real-valued vector $u$ and an integer-valued vector $z$,
according to their admissible values for $\XX$.
Furthermore, we consider the following blanket assumptions.
\begin{mybox}
	\begin{assumption}\label[assumption]{ass:P}
		With regard to \eqref{eq:P},
		\begin{itemize}
			\item functions $f$ and $c$ are continuously differentiable with locally Lipschitz derivatives;
			\item admissible values for the integer-valued decision variables in $z$ lie in a bounded set.
		\end{itemize}
	\end{assumption}
\end{mybox}

This work proposes a solution technique for (discretized) hybrid OCPs based on a safeguarded augmented Lagrangian (AL) approach, traditionally confined to nonlinear programming.
Adopting the AL framework of \cite{birgin2014practical}, the constrained problem \eqref{eq:P} is solved through a sequence of simpler subproblems.  
This methodology is motivated by sequential (partially) unconstrained minimization schemes, including penalty and barrier methods; see \cite{fiacco1968nonlinear}.
The AL subproblems are tackled by the tailored solver proposed in \cite{demarchi2024mixed}.
Such a practical AL algorithm for MINLPs is our main contribution, intended to seek good local, possibly non-global, solutions to \eqref{eq:P}.

The remainder of this work is structured as follows.
A summary of the optimality concept is described in \cref{sec:optimalityConditions}, while algorithmic details are provided in \cref{sec:algorithmicFramework}.
Hybrid OCPs are addressed with our method in \cref{sec:numericalResults}, where numerical results demonstrate its modelling versatility and efficacy.

\section{Optimality Conditions}
\label{sec:optimalityConditions}
This section is dedicated to developing suitable solution concepts for \eqref{eq:P}.
First of all, we seek \emph{feasible} points, namely some $\bar{x} \in \XX$ such that $c(\bar{x}) \in \CC$.
Then, a \emph{global} minimizer $\bar{x}$ for \eqref{eq:P} can be readily characterized by
\begin{align*}
	&\bar{x} \in \XX ,\quad
	c(\bar{x}) \in \CC , \quad
	\forall x\in\XX ,\, c(x)\in\CC \colon~ f(\bar{x}) \leq f(x) .
\end{align*}
Since we aim at developing an affordable numerical method, possibly at the price of global optimality,
we are interested in \emph{local} optimality notions too.
However, as these are sensitive to how neighborhoods are defined,
local notions can be weak and fragile in the mixed-integer setting considered here.
Indeed, considering simple balls, sufficiently small neighborhoods of feasible 
points may contain only one integer configuration,
making it locally optimal, in a weak sense.
To avoid this naive declaration of local optimality,
we develop some stronger conditions by
constructing neighborhoods based on
a polyhedral norm applied to the real-valued components only,
denoted $\normlp{\cdot}$ as in
\cite[\S 2]{demarchi2024mixed}.
For instance, it can take on the forms
\begin{align*}
	\normlp{x} \coloneqq{}& \sum\{ |x_i| \,|\, i\notin\II \} &\text{for }&\ell_1\text{-type}, \\
	\text{or}\quad
	\normlp{x} \coloneqq{}& \max\{ |x_i| \,|\, i\notin\II \} &\text{for }&\ell_\infty\text{-type} .
\end{align*}
Although not a norm, $\normlp{\cdot}$ induces compact neighborhoods according to
\[
\lpball(\bar{x},\Delta) \coloneqq \{ x\in\spaceX \,|\, \normlp{x-\bar{x}}\leq\Delta \},
\]
owing to the boundedness of integer-valued variables by \cref{ass:P}.
Now, analogously to the ``unconstrained'' case,
local minimizers for \eqref{eq:P} can be defined as follows.
\begin{mybox}
	\begin{definition}\label[definition]{def:minimizer}
		A point $\bar{x}\in\spaceX$ is called a \emphdef{local minimizer}
		for \eqref{eq:P} if it is feasible
		and
		there exists $\Delta>0$ such that
		$f(\bar{x}) \leq f(x)$ for all feasible $x\in\lpball(\bar{x},\Delta)$.
		If the latter property additionally holds for all $\Delta>0$, then
		$\bar{x}$ is called a \emphdef{global minimizer} for \eqref{eq:P}.
	\end{definition}
\end{mybox}

Notice that non-polyhedral norms for $\normlp{\cdot}$ are deliberately avoided so that, in practice, local models of \eqref{eq:P} take the form of mixed-integer \emph{linear} programs.
However, this specification is not necessary to the definition of suitable neighborhoods nor to identify stronger local minimizers.

We now focus on first-order necessary optimality conditions for \eqref{eq:P}.
These should be useful and practical in order to characterize and numerically detect
points that are (at least) candidates for global (or local) minimizers \cite[Chapter 3]{birgin2014practical}.

\subsection{Mixed-integer polyhedral constraints}
Let us first consider the minimization of some smooth function $\func{\lf}{\XX}{\R}$ over $\XX$:
\begin{align}
	\label{eq:Punc}
	\minimize\quad&\lf(x)&
	\wrt\quad&x\in\XX .
\end{align}
Recalling \cite[Def. 2.3, Prop. 2.4]{demarchi2024mixed}, a first-order optimality measure associated to \eqref{eq:Punc}
is defined
for all $\bar{x}\in\XX$ and $\Delta \geq 0$ by
\begin{equation}
	\label{eq:psimeas}
	\psimeas_{\lf,\XX}(\bar{x}, \Delta)
	\coloneqq
	\max_{\substack{x\in\XX\\\normlp{x - \bar{x}} \leq \Delta}} \innprod{\nabla \lf(\bar{x})}{\bar{x} - x}
	\geq
	0
	.
\end{equation}
Then, a first-order optimality concept for \eqref{eq:Punc} is that in the following definition,
which also provides an approximate counterpart thereof.
\begin{mybox}
	\begin{definition}\label[definition]{def:criticality}
		Given some $\varepsilon \geq 0$,
		a point $\bar{x} \in \XX$ is said \emphdef{$\varepsilon$-critical} for \eqref{eq:Punc} if there exists some $\Delta>0$ such that 
		$\psimeas_{\lf,\XX}(\bar{x}, \Delta) \leq \varepsilon$.
		A $0$-critical point is simply called \emphdef{critical}.
	\end{definition}
\end{mybox}
As discussed in \cite[\S 2.2]{demarchi2024mixed},
criticality is a necessary optimality condition for \eqref{eq:Punc} and, as such, plays a crucial role when it comes to designing an affordable numerical method for solving \eqref{eq:Punc}.
The denomination \emph{criticality}, in contrast to \emph{stationarity}, emphasizes that in general the former notion implies (hence, is stronger than) the latter.

\subsection{Nonlinear constraints}

The ``unconstrained'' notion of \cref{def:criticality} is useful to characterize solutions of a problem with simple constraints.
But what is a ``critical point'' for a problem such as
\eqref{eq:P}?
A stationarity characterization that resembles, at least in spirit, the classical Karush-Kuhn-Tucker (KKT) conditions in nonlinear programming
has been proposed and analysed by \cite[\S 2.2]{demarchi2024affordable}.

Let the Lagrangian function $\func{\Lagr}{\XX\times\spaceY}{\R}$ associated to \eqref{eq:P} be defined by
\begin{equation}
	\label{eq:lagrangian}
	\Lagr(x,y)
	\coloneqq
	f(x) + \innprod{y}{c(x)}
\end{equation}
where $y$ denotes a Lagrange multiplier.
From the viewpoint of nonlinear programming, where local minima correspond to saddle points of the Lagrangian functions, we consider the following notion for KKT-like points of \eqref{eq:P},
building upon \cref{def:criticality}.
\begin{mybox}
\begin{definition}\label[definition]{def:KKTcritical}
	A point $\bar{x} \in \spaceX$ is said \emphdef{KKT-critical} for \eqref{eq:P} if
	$\bar{x}\in\XX$ and there exists a multiplier $y \in \spaceY$ and some $\Delta>0$ such that
	\begin{equation*}
		\psimeas_{\Lagr(\cdot,y),\XX}(\bar{x},\Delta) = 0
		\qquad\text{and}\qquad
		y \in \normalcone_\CC(c(\bar{x})) .
	\end{equation*}
\end{definition}
\end{mybox}

Notice that KKT criticality demands feasibility since the normal cone $\normalcone_\CC(c(\bar{x}))$ to $\CC$ at $c(\bar{x})$ would be empty otherwise.
By \eqref{eq:psimeas}--\eqref{eq:lagrangian}, the first condition in \cref{def:KKTcritical}
can be rewritten as
\begin{equation*}
	\min_{\substack{x\in\XX\\\normlp{x - \bar{x}} \leq \Delta}} \innprod{\nabla f(\bar{x}) + c^\prime(\bar{x})^\top y}{x - \bar{x}}
	=
	0
	,
\end{equation*}
meaning that the Lagrangian function at $\bar{x}$ cannot be (locally) further minimized with respect to $x$ while maintaining MIL feasibility, in the sense of \cref{def:criticality}.
Thanks to the $\normlp{\cdot}$-localization, this condition is in general stronger than the inclusion
$0 \in \nabla_x \Lagr(\bar{x},y) + \limnormalcone_{\XX}(\bar{x})$,
corresponding to mere stationarity,
where $\limnormalcone_{\XX}(\bar{x})$ denotes the limiting normal cone to $\XX$ at $\bar{x}$.
Equivalence holds when the set $\XX$ is convex, which is not the case as soon as $\II\neq\emptyset$.

A criterion for numerical termination arises from
relaxing criticality and feasibility requirements
in \cref{def:KKTcritical}.
Analogously to $\varepsilon$-criticality of \cref{def:criticality},
we consider the following concept,
akin to \cite[Def. 3.2]{demarchi2024implicit}.
\begin{mybox}
	\begin{definition}\label[definition]{def:epsKKTcritical}
	Given some $\varepsilon \geq 0$,
	a point $\bar{x} \in \spaceX$ is said \emphdef{$\varepsilon$-KKT-critical} for \eqref{eq:P} if $\bar{x} \in \XX$ and there exists a multiplier $y \in \spaceY$, a vector $z\in\CC$, and some $\Delta>0$ such that
	\begin{align*}
		\psimeas_{\Lagr(\cdot,y),\XX}(\bar{x},\Delta) \leq{}& \varepsilon ,&
		y \in{}& \normalcone_\CC( z ) ,&
		\| c(\bar{x}) - z \| \leq{}& \varepsilon .
	\end{align*}
	\end{definition}
\end{mybox}
Consistently with \cref{def:KKTcritical}, a $0$-KKT-critical point is KKT-critical.
Possibly under some constraint qualifications (e.g., LICQ), KKT-criticality provides a necessary optimality condition for \eqref{eq:P}.
However, regardless of additional regularity, each local minimizer of \eqref{eq:P} is \emph{asymptotically} KKT-critical, cf. \cite[Thm 2.6]{demarchi2024affordable}.

\section{Algorithmic Framework}
\label{sec:algorithmicFramework}
We employ an augmented Lagrangian (AL) scheme to seek a numerical solution for \eqref{eq:P} under \cref{ass:P}.
This design choice is to exploit the fact that the AL framework maintains the original problem's structure and it does not rely on smoothness of the objective nor on convexity of the feasible set \cite{birgin2014practical}.
Building upon a sequence of shifted-penalty subproblems,
AL algorithms can be described with a nested structure.
An outer loop generates a sequence of primal-dual estimates,
adapts the penalty parameter,
and monitors convergence toward KKT-critical points.
Each subproblem associated to a set of parameters is addressed by an iterative procedure, effectively an inner loop.
Although the AL framework comprises several algorithms and variants,
all these share a core element:
the successive parametric minimization of an AL function.
The safeguarded AL scheme for \eqref{eq:P} stated in \cref{alg:ALM}
patterns that of \cite[\S 4.1]{birgin2014practical} and
proceeds with a sequence of subproblems of the form
\begin{align}
	\label{eq:ALsubproblem}
	\minimize\quad&\Lagr_\mu(x,\widehat{y})
	&
	\wrt\quad&x \in\XX
\end{align}
for some given penalty parameter $\mu>0$ and multiplier estimate $\widehat{y}\in\spaceY$.
Herein, the MIL constraints $x\in\XX$ are kept implicit, and the AL function is defined by
\begin{align}
	\Lagr_\mu(x,\widehat{y})
	\coloneqq{}&
	f(x) + \frac{1}{2 \mu} \dist_{\CC}^2\left( c(x) + \mu \widehat{y} \right) - \frac{\mu}{2} \|\widehat{y}\|^2 ,
\end{align}
cf. \cite[\S 4]{demarchi2024implicit}.
Feasibility of \eqref{eq:ALsubproblem} follows from $\XX$ being nonempty
whereas well-posedness is due to continuity of $\Lagr_\mu(\cdot,\widehat{y})$ and can be guaranteed, e.g., by coercivity or (level) boundedness arguments.
Notice that $\Lagr$ and $\Lagr_\mu$ are differentiable with respect to both primal and dual variables, thanks to convexity of $\CC$, with derivatives
\begin{subequations}
	\begin{align}
		\nabla_x \Lagr_\mu(x,\widehat{y})
		={}&
		\nabla f(x) + c^\prime(x)^\top y_\mu(x,\widehat{y}) , \\
		\nabla_y \Lagr_\mu(x,\widehat{y})
		={}&
		c(x) - s_\mu(x,\widehat{y})
	\end{align}
\end{subequations}
where
\begin{subequations}
	\begin{align}
		s_\mu(x,\widehat{y}) \coloneqq{}& \proj_{\CC}(c(x) + \mu \widehat{y}) , \\
		y_\mu(x,\widehat{y}) \coloneqq{}& \widehat{y} + \frac{c(x) - s_\mu(x,\widehat{y})}{\mu}
		.
	\end{align}
\end{subequations}
Furthermore, note that the input requirement in \cref{alg:ALM} is not a hard restriction:
given some $x^0 \notin \XX$, it is enough to project it onto $\XX$ by solving a MILP, for instance $\min_{x\in\XX} \|x-x^0\|_1$.

\begin{algorithm2e}[tb]
	\DontPrintSemicolon
	\KwIn{$x^0\in\XX$, $y^0\in\spaceY$, $\epsilon_{\textup{p}},\epsilon_{\textup{d}} > 0$,
		$\mu_1, \varepsilon_1 > 0$, $\kappa_\mu, \theta_\mu, \kappa_\varepsilon \in (0,1)$, $\Ybounded \subset \spaceY$ compact}
	\For{$j = 1,2\ldots$}{
		$\widehat{y}^j \gets \proj_{\Ybounded}(y^{j-1})$\label{step:ALM:ysafe}\;
		$x^j \gets$ compute an $\varepsilon_j$-critical point for $\Lagr_{\mu_j}(\cdot,\widehat{y}^j)$ over $\XX$ starting from $x^{j-1}$\label{step:ALM:subproblem}\;
		$s^j \gets \proj_{\CC}(c(x^j) + \mu_j \widehat{y}^j)$
		and
		$v^j \gets c(x^j) - s^j$\label{step:ALM:cviol}\;
		$y^j \gets \widehat{y}^j + v^j/\mu_j$ \label{step:ALM:y}\;
		\lIf{$\varepsilon_j \leq \epsilon_{\textup{d}}$ \KwAnd $\|v^j\| \leq \epsilon_{\textup{p}}$}{%
			\KwReturn $(x^j,y^j)$\label{step:ALM:termination}
		}
		\lIf{$j=1$ \KwOr $\|v^j\| \leq \max\{ \epsilon_{\textup{p}}, \theta_\mu \|v^{j-1}\|\}$}{%
			$\mu_{j+1} \gets \mu_j$, \KwElse $\mu_{j+1} \gets \kappa_\mu \mu_j$
		}\label{step:ALM:penaltyp}
		$\varepsilon_{j+1} \gets \max\{\epsilon_{\textup{d}}, \kappa_\varepsilon\varepsilon_j\}$\label{step:ALM:innertol}\;
	}
	\caption{Safeguarded augmented Lagrangian routine for \eqref{eq:P}}
	\label{alg:ALM}
\end{algorithm2e}

The major computational toll is taken by \cref{step:ALM:subproblem}, which aims to minimize $\Lagr_{\mu_j}(\cdot,\widehat{y}^j)$ over $\XX$.
In fact, only approximate criticality is required, and the sequential mixed-integer linearization scheme of \cite[Alg. 3.1]{demarchi2024mixed} can be readily adopted to compute a suitable $\varepsilon_j$-critical point $x^j$, in the sense of \cref{def:criticality}.
Although the performance of the algorithm might be improved by warm-starting the inner solver and guaranteeing descent behaviour,
these properties are not strictly required in the convergence analysis \cite{demarchi2024affordable}.
The classical first-order dual estimate update takes place at \cref{step:ALM:cviol,step:ALM:y}.
Given a (possibly inexact, first-order) solution $x$ to \eqref{eq:ALsubproblem}, the dual update rule in \cref{step:ALM:y}
is designed toward the identity
\begin{equation*}
	\nabla_x \Lagr_\mu(x,\widehat{y})
	={}
	\nabla f(x) + c^\prime(x)^\top y
	={}
	\nabla_x \Lagr(x,y) ,
\end{equation*}
as usual in AL methods.
This relationship allows to monitor the ``outer'' convergence for \eqref{eq:P} with the ``inner'' subproblem tolerance for \eqref{eq:ALsubproblem}.
Since the inclusion $y\in\normalcone_\CC(s)$ in \cref{def:epsKKTcritical} is always satisfied by construction of $s^j$ and $y^j$ in \cref{step:ALM:y,step:ALM:cviol},
approximate KKT criticality can be directly detected in \cref{step:ALM:termination} returning a pair $(x^j,y^j)$.
The dual safeguard takes place in \cref{step:ALM:ysafe},
where the compact set $\Ybounded$ can be a generic hyperbox or tailored to the constraint set $\CC$ to exploit additional known structures.
Finally, penalty parameter and inner tolerance are adapted in \cref{step:ALM:penaltyp,step:ALM:innertol} following classical update rules \cite[\S 12.4]{birgin2014practical}.
Initial values for steering primal-dual tolerance sequences can be user-specified or adaptively selected based on the initial infeasibility and criticality measure.
We shall mention also that, when a feasible initial point is available,
reset schemes are applicable and provide asymptotic feasibility guarantees
hardly attainable otherwise.

In analogy to magical steps \cite[\S 8.2]{birgin2014practical}, \cite[\S 4.1]{demarchi2024implicit} and feasibility pumps \cite{dambrosio2012storm}, additional steps can be optionally included in \cref{alg:ALM} to ease or accelerate the solution process.
Such routines could pursue feasibility (switching to the minimization of an infeasibility measure) or refine the incumbent solution by solving a nonlinear program with fixed integer-valued variables \cite[Rem. 3.1]{demarchi2024mixed}, \cite[\S 5]{kirches2022sequential}.

\begin{remark}
	A favorable attribute of AL methods is that if the AL subproblems are solved to approximate global optimality, then the iterates accumulate at feasible, globally optimal points \cite[Thm 5.2]{birgin2014practical}.
	Thus, whenever additional knowledge or problem structure allow to solve the AL subproblems to approximate global optimality, this general property of AL methods automatically provides global guarantees for \eqref{eq:P}.
\end{remark}

\section{Numerical Results}
\label{sec:numericalResults}
In this section, we present and compare the results achieved by the proposed numerical approach for two hybrid optimal control problems.
Before considering the simulations, we highlight some aspects regarding the computational procedure behind \cref{alg:ALM} (ALM).

ALM is initialized with zero dual estimate $y^0$, penalty parameter $\mu_1 = 10^{-2}$ and inner tolerance $\varepsilon_1 = 10^{-2}$.
The choice of $x^0$ is based on a more sophisticated approach.
Admissible points for the set $\XX$ are computed as projections onto it, in the $\ell_1$-norm sense.
Another routine was implemented to generate reasonable starting points and to refine solutions provided by ALM.
While integer variables are relaxed in the former case,
they are fixed (to an admissible value) in the latter case so that only real-valued variables are further optimized.
In both scenarios, one ends up with a nonlinear problem, which is then tackled by Ipopt \cite{waechter2006implementation}.
We use the numerical routine of \cite[\S 4]{demarchi2024mixed} for executing \cref{step:ALM:subproblem}, with Gurobi for solving the MIL subproblems there.

\subsection{Car with hysteretic turbo charger and drag}
The first problem handled by the proposed method is inspired by one from
\cite[\S 4.1]{demarchi2024mixed}.
The aim is to perform point-to-point one-dimensional motion planning of a car with hybrid nonlinear dynamics,
extending the previous model to accommodate a drag force term.
The underlying double-integrator point-mass model is equipped with a hysteretic turbo accelerator.
The car is described by its time-dependent position $x(t)$, velocity $v(t)$ and a turbo state $w(t) \in\{0,1\}$ for each time point $t \in [0, T]$, $T > 0$.
The control variables correspond to the inputs to the acceleration $a(t)$ and brake $b(t)$ pedals.
The turbo mode becomes active or inactive once the velocity exceeds $v^+ \coloneqq 10$ or falls below the threshold $v^- \coloneqq 5$, respectively.
An activated turbo mode entails an effective increase of the nominal thrust by three times but leaves the braking force unaffected.
These dynamics read
\begin{align*}
	\dot{x}(t) ={}& v(t) , &
	\dot{v}(t) ={}& \tau(w(t),a(t)) - b(t) - c_{\rm d} v(t)^2
\end{align*}
where $c_{\rm d} \geq 0$ denotes the drag coefficient, and
the thrust $\tau$ has two modes of operation depending on the turbo state:
$\tau(w,a) = a$ if $w = 0$, and $\tau(w,a) = 3a$ if $w = 1$.
The velocity is limited by $|v(t)|\leq v_{\max}\coloneqq 25$.
Boundary conditions, state and control bounds, and the minimum-control objective function are chosen following \cite[\S 4.1]{demarchi2024mixed}.
The same holds for the hysteresis characteristic describing the turbo behaviour,
whose logical implications are encoded as MIL constraints.

\paragraph*{Simulations}

Considering the final time $T = 10$, the problem is fully discretized using the explicit Euler scheme on a uniform time grid with $N=100$ intervals.
The resulting problem in the form \eqref{eq:P} has $100$ binary-valued optimization variables and $101$ nonlinear constraints.

\begin{figure}[tbh]
	\centering
	\includegraphics{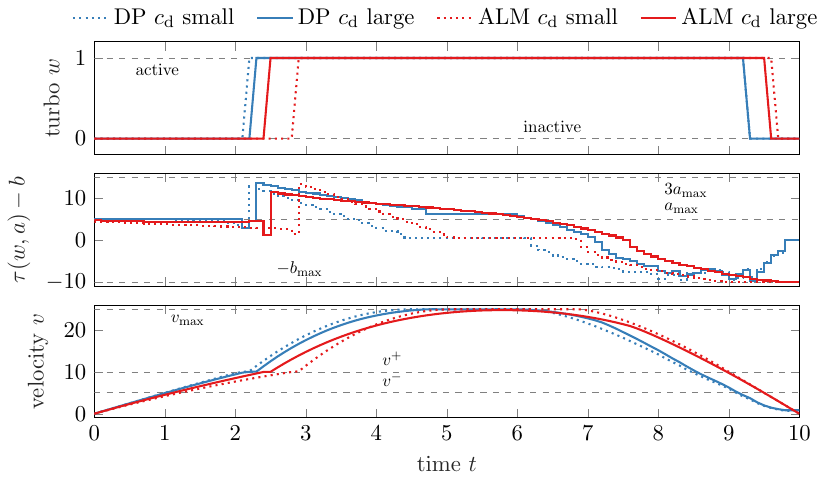}
	\caption{Solutions for the turbo car problem with different values of drag coefficient $c_{\rm d}\in\{10^{-3},10^{-2}\}$. Comparison of \cref{alg:ALM} (ALM), starting from an all-zero initial guess, and dynamic programming (DP).}
	\label{fig:turbo_car_drag}
\end{figure}

ALM is initialized with an all-zero primal input.
The solutions returned for different values of the drag coefficient $c_{\rm d} \in \{10^{-3},10^{-2}\}$ are displayed in \cref{fig:turbo_car_drag},
together with the solution obtained with dynamic programming (DP).
For DP, we used uniform discretization grids with 100 values for each (real-valued) state and control.
In all tests, boundary conditions and path constraints are approximately satisfied.
As expected, all solutions resemble the one presented in \cite{demarchi2024mixed} for the linear case with $c_{\rm d} = 0$, whereas more acceleration and less braking are demanded with stronger drag.
The DP solution exhibits some chattering behaviour in the braking phase---this result may be an artifact due to discretization.
However, finer discretizations appeared impractical in terms of runtime (on a standard laptop):
for each call, DP took minutes to compute, whereas ALM was in the order of seconds.

\begin{figure}[tbh]
	\centering
	\includegraphics{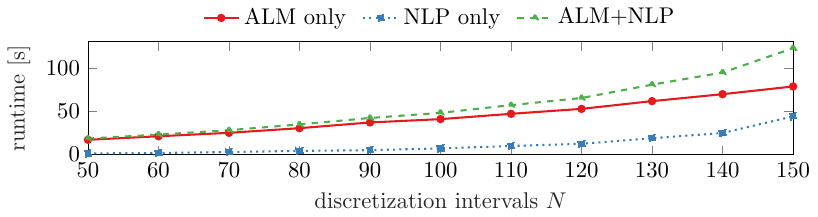}
	\caption{Runtimes for the turbo car problem with drag $c_{\rm d}=10^{-3}$ for different number $N$ of discretization intervals. Solutions obtained with ALM and subsequent fixed-integer NLP refinement.}
	\label{fig:turbo_car_runtime}
\end{figure}

For $c_{\rm d} = 10^{-3}$, the turbo car problem was in addition solved by ALM combined with a subsequent fixed-integer NLP refinement to examine the impact of the problem size on the computational effort.
As illustrated in \cref{fig:turbo_car_runtime}, the NLP runtime grows faster but is lower than ALM's, which however increases only linearly with $N$.

\subsection{Lotka-Volterra fishing with total variation}

Let us now consider the classical Lotka-Volterra fishing problem complemented with a total variation ($\totalvariation$) term \cite{sager2012mintoc}.
For two states $x_1$ and $x_2$, which can be interpreted as the biomasses of prey and predator populations, respectively, a binary optimal control is sought.
The dynamics are described with the model
\begin{align*}
	\dot{x}_1(t) ={}& x_1(t) - x_1(t) x_2(t) - c_1^\top w(t) , \\
	\dot{x}_2(t) ={}& x_1(t) x_2(t) - x_2(t) - c_2^\top w(t)
\end{align*}
for some given parameters $c_1, c_2\in\R^5$ and control inputs $w(t) \in \{0,1\}^5$.
For every time $t\in\R$, the pointwise constraint $\sum_{i=1}^5 w_i(t) = 1$ encodes the choice of one among five control options.
The primary objective is formulated as 
\begin{equation*}
	J(x)
	\coloneqq
	\min \quad \int_{0}^{T} \| x(t) - x_r \|^2 \mathrm{d}t
\end{equation*}
for tracking the reference value $x_r = (1,1)$, starting from the initial state $x_1(0) = 0.5$, $x_2(0)=0.7$ and with final time $T=12$ \cite[\S 5]{sager2011combinatorial}.

To limit the well-known chattering behaviour of the optimal control,
we monitor the total variation (in fact, its discrete counterpart, after time discretization)
\begin{equation*}
	\totalvariation(w) \coloneqq \frac{1}{2} \sum_{i=1}^5 \sum_{k=1}^{N-1} |w_{i,k+1} - w_{i,k}|
	,
\end{equation*}
possibly imposing an upper bound $\totalvariation(w) \leq U_{\totalvariation}$ or including a penalty term $\alpha_{\totalvariation} \totalvariation(w)$ in the objective.
By including additional $(N-1)$ real-valued variables and $4(N-1)$ linear inequality constraints for each binary-valued control, the nonsmooth $\totalvariation$ term can be expressed as a linear function of the unknowns.

\paragraph*{Simulations}

The mixed-integer OCP is fully discretized using the explicit Euler scheme over a uniform time grid with $N = 50$ intervals.
The resulting discretized problem contains $250$ integer-valued optimization variables and $51$ nonlinear constraints.

\begin{figure}[tbh]
	\centering
	\includegraphics{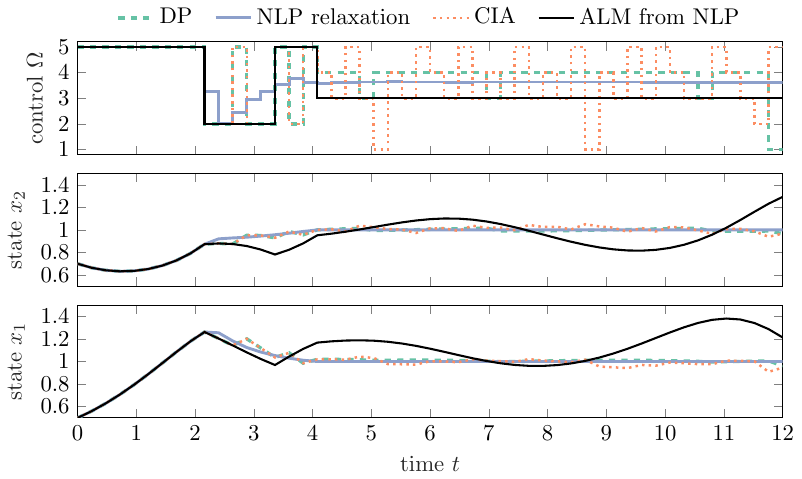}
	\caption{Solutions for the fishing problem without $\totalvariation$, starting CIA and ALM from the NLP relaxation.}
	\label{fig:fishing_tv}
\end{figure}

The initial guess values are set to $1$ for all variables and the problem (with no TV cost nor limit) is first considered with relaxed integrality and solved by Ipopt.
Such solution is then used as input for the combinatorial integral approximation (CIA) \cite{sager2011combinatorial}, computed via \texttt{pycombina} \cite{buerger2020pycombina}, and for warm-starting ALM.
The refined solutions are displayed in \cref{fig:fishing_tv} and summarized in \cref{tab:fishing_tv}, together with the globally optimal solution obtained with DP (on a uniform grid with 100 values for each state).
We report the effective control $\Omega \coloneqq \sum_{i=1}^5 i w_i$, instead of the vector-valued $w$, to better display the selected fishing option.
When warm-started from the NLP relaxation, the infeasible initial guess is first projected onto $\XX$ before running ALM; the returned solution corresponds to such projection and does not provide good tracking performance.
Similarly, when warm-started from the DP or CIA solution, ALM returns the same control sequence, validating their local optimality.

\begin{table}[tbh]
	\renewcommand{\arraystretch}{1.15}
	\centering
	\caption{Results for the fishing problem with and without $\totalvariation$: upper bound $U_{\totalvariation}$ on the total variation $\totalvariation$, tracking cost $J$. When not indicated, ALM starts from the CIA solution.}
	\label{tab:fishing_tv}
	\begin{tabular}{c|ccc}
		\hline
		solver & $U_{\totalvariation}$ & $\totalvariation(w)$ & $J(x)$ \\
		\hline
		DP & $\infty$ & 14 & 0.388 \\
		ALM from DP & $\infty$& 14 & 0.388 \\
		NLP relaxation & $\infty$ & 2.6 & 0.381 \\
		ALM from NLP & $\infty$ & 3 & 0.834 \\
		CIA from NLP & $\infty$ & 38 & 0.398 \\
		ALM & $\infty$ & 38 & 0.398 \\
		\hline
	\end{tabular}
	\quad
	\begin{tabular}{c|ccc}
		\hline
		solver & $U_{\totalvariation}$ & $\totalvariation(w)$ & $J(x)$ \\
		\hline
		ALM & 38 & 38 & 0.398 \\
		ALM & 33 & 33 & 0.415 \\
		ALM & 28 & 28 & 0.825 \\
		ALM & 23 & 23 & 0.549 \\
		ALM & 18 & 18 & 1.206 \\
		ALM & 13 & 13 & 0.601 \\
		\hline
	\end{tabular}
\end{table}

\begin{figure}[tbh]
	\centering
	\includegraphics{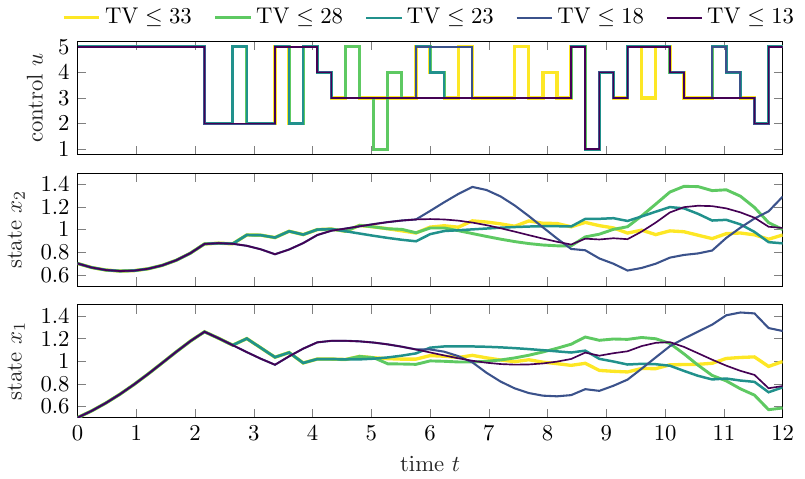}
	\caption{Solutions for the fishing problem with limited TV, starting ALM from the CIA solution.}
	\label{fig:fishing_tv_ALM}
\end{figure}

We now consider the performance of ALM when limiting the control $\totalvariation$.
Always starting from the CIA solution, ALM is executed with decreasing values $U_{\totalvariation}$ of maximum admissible $\totalvariation$.
The corresponding results are shown in \cref{fig:fishing_tv_ALM} and collected in \cref{tab:fishing_tv}.
Note that since the ALM solution coincides with the CIA solution for $U_{\totalvariation} = 38$, this case is omitted in the plots from \cref{fig:fishing_tv_ALM}.

Expectedly, limiting the (discretized) $\totalvariation$ helps avoiding a chattering behaviour of the optimal solution but at the cost of impairing the tracking performance.
In all cases, the returned solution had the maximum TV allowed.
Moreover, \cref{tab:fishing_tv} reports that the solution for $U_{\totalvariation}=23$ (resp. $U_{\totalvariation}=13$) had a better tracking objective than the one for $U_{\totalvariation}=28$ ($U_{\totalvariation}=18$), testifying the local nature of ALM when coupled with the local inner solver of \cite{demarchi2024mixed}.

\section{Conclusions}

This work presents an efficient method to address hybrid optimal control problems without any integrality relaxations.
The proposed approach is based on the augmented Lagrangian framework combined with a routine for nonlinear optimization with mixed-integer linear constraints.
The robustness and performance of the algorithm were assessed on two problem examples from different fields of application, validating the proposed strategy.
Future research should consider the integration of other techniques for handling nonlinear constraints.

\subsection*{Acknowledgements}
\TheAcknowledgements

%---------- References
\phantomsection
\addcontentsline{toc}{section}{References}%
\small
\bibliographystyle{plain}
\bibliography{biblio}

\end{document}